\documentclass[12pt]{article}

\usepackage[paper=a4paper,left=24mm,right=24mm,top=24mm,bottom=34mm]{geometry}
\usepackage{amsfonts,amsmath,amssymb}
\usepackage{bbm,booktabs,graphicx,multirow,setspace}
\usepackage{color}
\usepackage[hyphens]{url}
\usepackage[english]{babel}
\usepackage[utf8]{inputenc}
\usepackage[round]{natbib}

\newcommand{\pr}{\text{pr}}

\begin{document}

\title{When is the mode functional the Bayes classifier?}

\author{Tilmann Gneiting
        \\
        {\small Heidelberg Institute for Theoretical Studies and Karlsruhe Institute of Technology} 
        }

\date{\small \today}

\maketitle

\begin{abstract}
In classification problems, the mode of the conditional probability
distribution, i.e., the most probable category, is the Bayes
classifier under zero-one or misclassification loss.  Under any other
cost structure, the mode fails to persist.
\end{abstract}

Consider a finite number of categories or classes labeled $1, \ldots,
k$.  Let the random variable $Y$ denote the class label, and let $X$
be a random covariate or feature vector, for a unit at hand. A
{\em probabilistic classifier}\/ is a conditional probability distribution,
\[
p(i \mid x) = \pr(Y = i \mid X = x) \quad \text{for} \quad i = 1, \ldots, k, 
\]
which in typical practice is estimated from training data.  In
contrast, a deterministic classifier or {\em decision rule}\/ $G(x)$ assigns a
single class label to any realized feature.

The most common way of converting a probabilistic classifier into a
decision rule is to use the {\em mode}\/ functional $\bar{G}$, which assigns 
the most probable class label, i.e., 
\begin{equation}  \label{eq:mode}
\bar{G}(x) = i \quad \mbox{if} \quad p(i \mid x) = \max_{i' = 1, \ldots, k} \, p(i' \mid x),
\end{equation}
where ties are resolved by randomization. If $L(i,j)$ denotes the loss
or cost when $G(x) = i$ and class $j$ realizes, where $i, j = 1, \ldots,
k$, the associated {\em Bayes classifier}\/ or optimal decision rule
$\hat{G}$ assigns the class that minimizes the expected loss, i.e.,
\begin{equation}  \label{eq:Bayes}
\hat{G}(x) = i \quad \mbox{if} \quad
\sum_{j=1}^k  L(i,j) \, p(j \mid x) = \min_{i' = 1, \ldots, k} \, \sum_{j=1}^k  L(i',j) \, p(j \mid x). 
\end{equation}
The literature typically studies classification problems under
zero-one or misclassification loss, where $L(i,j) = 0$ if $i = j$ and
$L(i,j) = 1$ if $i \not= j$, and it is well known that the mode
functional is the associated Bayes classifier (Hastie et al.~2009, p.~21).

However, it has been argued that misclassification loss ``is rarely
what users of classification methods really want'' (Hand 1997, p.~7)
and that one would ``be hard pressed to find an application in which
the costs of different kinds of errors were the same'' (Witten et
al.~2011, p.~164).  Witten et al.~(2011, p.~167) further note that
under cost structures other than misclassification loss the Bayes
classifier ``might be different'' from the mode.  As we now show, it
will in fact be different, in the sense that under any other loss
structure the mode fails to persist.

To demonstrate this, we invoke the reasonableness condition of Elkan
(2001) and assume that $L(i,j) \geq L(i,i)$ for $i, j = 1, \ldots, k$, 
with at least one of the inequalities being strict.  Adding constants
columnwise concerns costs that depend on the outcome only, and
multiplying all entries of the loss matrix by a positive number merely
changes the monetary unit.  Therefore, we may restrict attention to
loss or cost matrices for which $L(i,i) = 0$, $L(i,j) \geq 0$ and
$\sum_{i \not= j} L(i,j) = k(k-1)$.  In other words, the diagonal
elements vanish, and the off-diagonal entries are nonnegative and
average to one.

In the binary case $k = 2$ we thus consider cost matrices of the
form
\[
\left( \begin{array}{cc} 0 & \; c \\ 2-c & \; 0 \end{array} \right) \! ,
\] 
where $0 \leq c \leq 2$.  The optimal decision is $\hat{G}(x) = 1$ if
$p(1 \mid x) \geq c/2$.  If $c < 1$ and $c/2 \leq p(1 \mid x) < 1/2$,
the mode fails to be optimal; if $c > 1$, an analogous argument
applies.  Hence the mode functional is the Bayes classifier under
zero-one loss only.

In the ternary case $k = 3$ we may restrict attention to symmetric
cost matrices of the form
\begin{equation}  \label{eq:ternary}
\left( \begin{array}{ccccc} 0 && a && b \\ a && 0 && 3 - a - b \\ b && 3 - a - b && 0 \end{array} \right) \! ,
\end{equation}
where $0 \leq a \leq 3$ and $0 \leq b \leq 3-a$, for if the cost
matrix is asymmetric, the above arguments apply to a principal
submatrix.  The optimal decision under the cost matrix
\eqref{eq:ternary} is $\hat{G}(x) = 1$ if
\[ 
2b \, p(1 \mid x) \geq (2a-3) \, p(2 \mid x) + b, \quad 
2a \, p(1 \mid x) \geq (2b-3) \: p(3 \mid x) + a;
\]
other cases are handled analogously. If $a = b = 1$, we recover
zero-one loss, and the inequalities reduce to the conditions for the
mode.  Else, they yield functionals other than the mode.

When $k \geq 4$ we see from the ternary case that a necessary
condition for the cost matrix to yield the mode functional as Bayes
classifier is that every $3 \times 3$ principal submatrix be of the
form 
\[
\left( \begin{array}{ccccc} 0 && c && c \\ c && 0 && c \\ c && c && 0 \end{array} \right) \! ,
\]
where $0 \leq c \leq k(k-1)/6$. Considering successive principal
submatrices, and iterating the argument, we see that in fact $c =
1$. Hence, the mode functional \eqref{eq:mode} is the Bayes classifier
\eqref{eq:Bayes} under zero-one loss only, subject to the above
assumptions. 

This result complements findings by Heinrich (2014) in the case of a
continuous outcome, where it is not possible to find a loss function
under which the mode functional is the Bayes predictor.  In the
discrete setting considered here, the rife failure of the most
probable value to minimize expected loss may urge practitioners to
work with probabilistic classifiers in lieu of deterministic decision
rules, as advocated powerfully by Harrell (2015, Section 1.3).

\subsection*{Acknowledgements}

This work was funded by the European Union Seventh Framework Programme
under grant agreement 290976.  The author is grateful for
infrastructural support by the Klaus Tschira Foundation and thanks
Werner Ehm, Alexander Jordan and Michael Strube for instructive
discussions.

\subsection*{References}

\newenvironment{reflist}{\begin{list}{}{\itemsep 0mm \parsep 1mm
\listparindent -7mm \leftmargin 7mm} \item \ }{\end{list}}

\vspace{-6.5mm}
\begin{reflist}

Elkan, C.~(2001). The foundations of cost-sensitive learning. In {\em
  Proceedings of the Seventeenth International Joint Conference on
  Artificial Intelligence}, Morgan Kaufmann, San Francisco,
pp.~973--978.

Hand, D.~J.~(1997). {\em Construction and Assessment of Classification
  Rules}.  Wiley, Chichester.

Harrell, F.~E.~(2015).  {\em Regression Modeling Strategies}, 2nd
edition.  Springer, Cham.

Hastie, T., Tibshirani, R.~and Friedman, J.~(2009).  {\em The Elements
  of Statistical Learning}, 2nd edition.  Springer, New York.

Heinrich, C.~(2014). The mode functional is not elicitable. Biometrika
{\bf 101}, 245--251.

Witten, I.~H., Frank, E.~and Hall, M.~A.~(2011).  {\em Data Mining},
3rd edition.  Elsevier Morgan Kaufmann, Amsterdam.

\end{reflist}

\end{document}